\documentclass[10pt]{article}
\usepackage[utf8]{inputenc}
\usepackage{latexsym,amssymb,amsmath,url,authblk}

\def\N{\mathbb{N}}
\def\Q{\mathbb{Q}}

\def\R{\mathbb{R}}

\def\proof{\par\noindent{\em Proof. }}
\def\eproof{\hfill{$\Box$}\bigskip}

\def\ds{\dots}
\def\sus{\subset}

\def\de{\delta}

\def\cc{\colon}
\def\ep{\varepsilon}

\newtheorem{thm}{Theorem}[section]
\newtheorem{prop}[thm]{Proposition}
\newtheorem{cor}[thm]{Corollary}

\newtheorem{defi}[thm]{Definition}

\title{Cardinalities in proofreading}
\author{M. Klazar\footnote{{\tt klazar@kam.mff.cuni.cz}}\\
KAM MFF UK Praha}
\date{\today}

\begin{document}

\maketitle
{\em Dedicated to the memory of my mother Blanka Klazarov\'a (1937--2026)}

\begin{abstract}
In 1975, G.~P\'olya suggested that if two proofreaders found $a$ and 
$b$ errors in 
a~text, of which $c$ errors were found by both of them, then a~reasonable approximation of 
the unknown number $e$ of all errors is $e\approx ab/c$. We
justify this formula by constructing a~realistic model of proofreaders and 
estimating the efficiency of this model with the help of the bound on 
large deviations in the Probabilistic Method. In conclusion 
we discuss the distinction between realistic and probabilistic models 
of problems. 
\end{abstract}

\section{Introduction. $(p,n)$-proofreaders}

In \cite{poly}, G.~P\'olya considered the following problem (see also \cite{alex,ande,ferg_hard,liu}). Suppose that a~proofreader found $a$ errors in 
a~text, that another proofreader found $b$ errors, and that the
number of errors detected by both of them is $c$. With these
data, how could one approximate the unknown number $e$ of all errors in the text? P\'olya's solution is this. We have $a\approx pe$ and $b\approx qe$, where $p$ and $q$ 
in $(0,1]$ are the probabilities of error detection for the two 
proofreaders. Their independence gives $c\approx pqe$. Thus we get the approximation
$$
e=\frac{peqe}{pqe}\approx\frac{ab}{c}\,.
$$
P\'olya's interesting argument is not supported by any rigorous 
theorem. In this article we provide one. It is 
Theorem~\ref{thm_physPrrea} which we prove in Section~\ref{sec_proof}. In 
Section~\ref{sec_concl} we discuss the two kinds of mathematical models 
one can devise to solve problems:  realistic models and probabilistic 
models. 

For $n$ in $\N=\{1,2,\ds\}$ we set $[n]=\{1,2,\ds,n\}$. Let $B_n=
\{0,1\}^n$ be the set of binary $n$-tuples 
$$
\overline{b}=(\overline{b}_1,\,\overline{b}_2,\,
\ds,\,\overline{b}_n),\ \overline{b}_i\in\{0,\,1\}\,.
$$
We denote by 
$${\textstyle
\|\overline{b}\|=\sum_{i=1}^n\overline{b}_i
}
$$ 
the number of $1$s in $\overline{b}$. 
For a~finite set $A$ we denote by $|A|$ ($\in\N_0=
\{0,1,\ds\}$) the number of its elements. For example, 
$|B_n|=2^n$. For a~tuple $\overline{s}$ we denote by $|\overline{s}|$
its length, the number of terms in $\overline{s}$. By $\Q$ we denote the set of fractions.

Let 
$p=k/l$ be in $(0,1]
\cap\Q$, with $k,l\in\N$, $k\le l$ and coprime $k,l$, and let $n\in\N$. 
We define the key discrete structure of our article. 

\begin{defi}[$(p,n)$-proofreaders]\label{def_pnPrRe}
A~$(p,n)$-proofreader $\mathcal{P}$ is the pair $\mathcal{P}=
(P,\overline{s})$ with the following components which are uniquely 
determined by $p$ and $n$.
\begin{enumerate}
\item $P=P(\overline{a},\overline{b})\cc B_n\times B_n\to B_n$ is a~map given by
$$
P(\overline{a},\,\overline{b})=
(\overline{a}_1\overline{b}_1,\,\overline{a}_2\overline{b}_2,\,\ds,\,\overline{a}_n\overline{b}_n)\,.
$$
\item The $N$-tuple of $n$-tuples
$$
\overline{s}=(\overline{s}_1,\,\overline{s}_2,\,
\ds,\,\overline{s}_N)\,, 
$$
where $\overline{s}_j\in B_n$ and $N=l^n$, is defined as follows. We 
begin with the map $f\cc[l]\to\{0,1\}$ given 
by $f(x)=1$ for $x\in[k]$ and $f(x)=0$ for $x\in[l]\setminus[k]$.
We extend it to the map $f\cc[l]^n\to B_n$ by
$$
f(a_1,\,a_2,\,\ds,\,a_n)=
(f(a_1),\,f(a_2),\,\ds,\,f(a_n))\,.
$$
We take the lexicographic ordering 
$\overline{u}_1$, $\overline{u}_2$, $\ds$, $\overline{u}_N$ of $[l]^n$, so that $\overline{u}_1=(1,1,\ds,1)$, $\overline{u}_2=
(1,1,\ds,1,2)$, $\ds$, $\overline{u}_N=(l,l,\ds,l)$. Finally, we set
$$
\overline{s}=(\overline{s}_1,\,\overline{s}_2,\,\ds,\,\overline{s}_N)=
(f(\overline{u}_1),\,f(\overline{u}_2),\,\ds,\,f(\overline{u}_N))\,.
$$
We call the $N$-tuple $\overline{s}$ the runs of $\mathcal{P}$.
\end{enumerate}
\end{defi}
Any
ordering of $[l]^n$ would do. We chose the lexicographic one 
just concreteness sake.

We show that the proofreaders introduced in Definition~\ref{def_pnPrRe} 
work according to the parameters $p$ and $n$. The parameter $p$ is the 
efficiency (probability of error detection) of the proofreader, and 
the parameter $n$ is the length of the  
text submitted for proofreading. $\mathcal{P}=(P,\overline{s})$ works in $N$ runs 
$\overline{s}_j$, $j\in[N]$. If $(\overline{s}_j)_i=1$ then in the $j$-th run  
$\mathcal{P}$ pays attention to the $i$-th position in 
the text, and if $(\overline{s}_j)_i=0$
then $\mathcal{P}$ ignores it. So if a~text $\overline{t}\in 
B_n$ (in which $1$s are errors) is submitted for proofreading, then $\mathcal{P}$ in the 
run $\overline{s}_j$ detects
$$
\|P(\overline{s}_j,\,\overline{t})\|\ \ (\in\N_0) 
$$
errors in $\overline{t}$. The 
next proposition shows that for every set $I\sus[n]$ of 
positions in the text, $\mathcal{P}$ pays attention to all positions in 
$I$ in exactly
$p^{|I|}\cdot N$ of the $N$ runs. 

\begin{prop}\label{prop_runs}
Let 
$\overline{s}=(\overline{s}_1,\overline{s}_2,\ds,\overline{s}_N)$
be the runs of the $(p,n)$-proofreader. Then for every set $I\sus[n]$ we have
$$
N^{-1}\cdot\big|
\big\{
j\in[N]\cc\;\text{$(\overline{ s}_j)_i=1$ for every $i\in I$}\}\big|=p^{|I|}\,.
$$
\end{prop}
\proof
Let $\overline{s}$, $p$, $n$, and $I$ be as stated. It is easy to see from the 
definition of $\overline{s}$ that the number of indices $j\in[N]$ such that $(\overline{s}_j)_i=1$
for every $i\in I$ equals
$$
k^{|I|}\cdot l^{n-|I|} 
$$
because it is the number of words of length $n$ over the alphabet $[l]$ that have at every position $i\in I$ a~letter from $[k]$ ($\sus[l]$) .
But then
$$
N^{-1}\cdot k^{|I|}\cdot l^{n-|I|}=
\frac{k^{|I|}\cdot l^{n-|I|}}{l^n}=\bigg(\frac{k}{l}\bigg)^{|I|}=p^{|I|}\,.
$$
\eproof

For P\'olya's argument
we need two independent proofreaders. Suppose that
$p,q\in(0,1]\cap\Q$ and that $n\in\N$. 

\begin{defi}[two $(p,q,n)$-proofreaders]\label{def_twoPrRe}
Two $(p,q,n)$-proofreaders is the pair $\mathcal{R}=(P,\overline{S})$ 
with the following components  which are uniquely 
determined by $p$, $q$ and $n$.
\begin{enumerate}
\item The map $P$ is the same as in Definition~\ref{def_pnPrRe}.
\item We construct the $MN$-tuple  $\overline{S}$ of pairs of $n$-tuples  
as follows. Let $(P,\overline{s})$
be the $(p,n)$-proofreader, $(P,\overline{u})$ be the 
$(q,n)$-proofreader, and let $M=|\overline{s}|$ and 
$N=|\overline{u}|$. We set 
\begin{eqnarray*}
\overline{S}&=&(\overline{S}_1,\,\overline{S}_2,\,\ds,\,\overline{S}_{MN})\\
&=&((\overline{s}_1,\,\overline{u}_1),\,(\overline{s}_1,\,\overline{u}_2),\,
\ds,\,(\overline{s}_1,\,\overline{u}_N),\,
(\overline{s}_2,\,\overline{u}_1),\,(\overline{s}_2,\,\overline{u}_2),\,
\ds,\,(\overline{s}_2,\,
\overline{u}_N)\,,\\&&\ds,\,\ds,\,
(\overline{s}_M,\,\overline{u}_1),\,(\overline{s}_M,\,\overline{u}_2),\,\ds,\,(\overline{s}_M,\,
\overline{u}_N))\ \ (\in(B_n\times B_n)^{MN})\,.
\end{eqnarray*}
We call the $MN$-tuple $\overline{S}$ the runs of $\mathcal{R}$.
\end{enumerate}
\end{defi}
Any ordering of the $MN$ pairs 
$(\overline{s}_j,\overline{u}_{j'})$ would do, but we again want to have 
a~uniquely defined object. For any $j\in[MN]$ we write
$$
\overline{S}_j=\big(\overline{s}_{v_j},\,\overline{u}_{w_j}\big)\,.
$$

In the
next proposition, we show that 
for every set $I\sus[n]$ of 
positions in the text submitted for proofreading, the two proofreaders 
pay attention to these positions mutually independently.

\begin{prop}\label{prop_twoProof}
Let
$\overline{S}=((\overline{s}_1,\overline{u}_1),
(\overline{s}_1,\overline{u}_2),
\ds,\ds,(\overline{s}_M,\,\overline{u}_N))$
be the runs of two $(p,q,n)$-proofreaders. The following holds.
\begin{enumerate}
\item For every set $I\sus[n]$ we have
$$
{\textstyle
\frac{1}{MN}\cdot\big|
\big\{
j\in[MN]\cc\;\text{$(\overline{s}_{v_j})_i=1$ for every $i\in I$}\big\}\big|=p^{|I|}\,.
}
$$
\item For every set $I\sus[n]$ we have
$$
{\textstyle
\frac{1}{MN}\cdot\big|
\big\{
j\in[MN]\cc\;\text{$(\overline{u}_{w_j})_i=1$  for every $i\in I$}\big\}\big|=q^{|I|}\,.
}
$$
\item For every set $I\sus[n]$ we have
$$
{\textstyle
\frac{1}{MN}\cdot\big|
\big\{
j\in[MN]\cc\;\text{$(\overline{s}_{v_j})_i=(\overline{u}_{w_j})_i=1$ for every $i\in I$}\big\}\big|
=(pq)^{|I|}\,.
}
$$
\end{enumerate}
\end{prop}
\proof
Let $\overline{S}$ and $I$ be as stated. 

1.~Let $r$ be the number of indices 
$j\in[M]$ such that $(\overline{s}_j)_i=1$ for every 
$i\in I$. By Proposition~\ref{prop_runs}, 
$\frac{r}{M}=p^{|I|}$. Thus the displayed product equals
$$
{\textstyle
\frac{1}{MN}\cdot rN=\frac{r}{M}=
p^{|I|}\,.
}
$$

2.~The argument is symmetric, with the number $s$ of indices
$j\in[N]$ such that $(\overline{u}_j)_i=1$ for every $i\in I$. 

3.~Now the displayed product equals
$$
{\textstyle
\frac{1}{MN}\cdot rs=\frac{r}{M}\cdot\frac{s}{N}=p^{|I|}\cdot
q^{|I|}=(pq)^{|I|}\,.
}
$$
\eproof

The following theorem is one of our main results.

\begin{thm}\label{thm_physPrrea}
Let $\mathcal{R}=(P,\overline{S})$, where 
$$
\overline{S}=((\overline{s}_1,\,\overline{u}_1),\,(\overline{s}_1,\,\overline{u}_2),\,
\ds,\,\ds,\,(\overline{s}_M,\,\overline{u}_N))\,,
$$ 
be two $(p,q,n)$-proofreaders and let
$$
{\textstyle
c(p,\,q)=\frac{1}{2}\big(\frac{1}{p}+\frac{1}{q}+\frac{2}{pq}\big)^{-1} \ \ (\in(0,\,\frac{1}{8}])\,.
}
$$
Then for every text $\overline{t}\in B_n$
with $\|\overline{t}\|=m\ge1$ errors and for every real number $a\in[0,m]$ we have
\begin{eqnarray*}
&&{\textstyle
\frac{1}{MN}}\cdot\big|\big
\{j\in[MN]\cc\;
\left|\frac{\|P(\overline{s}_{v_j},\,\overline{t})\|\cdot\|P(\overline{u}_{w_j},\,\overline{t})\|}{\|P(\overline{s}_{v_j},\,P(\overline{u}_{w_j},\,\overline{t}))\|}-m
\right|\le a\big\}\big|\ge\\
&&{\textstyle\ge1-
6\exp\big(-\frac{2c(p,\,q)^2a^2}{m}\big)\,,
}
\end{eqnarray*}
where for $\|P(\overline{s}_{v_j},P(\overline{u}_{w_j},\overline{t}))\|=0$ 
we define the fraction $\frac{\cdots}{0}$ as $+\infty$.
\end{thm}
The theorem says that the distance between P\'olya's approximation and the actual number 
$m$ of errors is at most 
$a$ for at least the displayed fraction $1-\ds$ of the $MN$ runs of 
$\mathcal{R}$.

The theorem has two corollaries. We
minimize the distance of the
approximation from the actual number of errors. Then we 
maximize the fraction of runs for which the distance is still 
reasonably small. The proofs are immediate, and we omit them.

\begin{cor}\label{cor_1stCor}
In the setup of Theorem~\ref{thm_physPrrea}, if $\omega\cc\N\to\N$ is a~function 
such that 
$$
\omega(m)\le\sqrt{m}\,\text{ and }\,\omega(m)\to\infty\,, 
$$
then for every text $\overline{t}\in B_n$
with $\|\overline{t}\|=m\ge1$ errors we have
\begin{eqnarray*}
&&{\textstyle
\frac{1}{MN}}\cdot\big|
\big\{j\in[MN]\cc\;\left|
\frac{\|P(\overline{s}_{v_j},\,\overline{t})\|\cdot\|P(\overline{u}_{w_j},\,\overline{t})\|}{\|P(\overline{s}_{v_j},\,P(\overline{u}_{w_j}
,\,\overline{t}))\|}-m
\right|\le\omega(m)\sqrt{m}\big\}\big|\ge\\
&&\ge1-6\exp(-2c(p,\,q)^2\omega(m)^2)\to1\ \ (m\to\infty)\,.
\end{eqnarray*}
\end{cor}

\begin{cor}\label{cor_2ndCor}
In the setup of Theorem~\ref{thm_physPrrea}, if $\ep\in(0,1)$ then for every text $\overline{t}\in B_n$
with $\|\overline{t}\|=m\ge1$ errors we have
\begin{eqnarray*}
&&{\textstyle
\frac{1}{MN}}\cdot\big|\big\{j\in[MN]\cc\;\left|
\frac{\|P(\overline{s}_{v_j},\,\overline{t})\|\cdot\|P(\overline{u}_{w_j},\,\overline{t})\|}{\|P(\overline{s}_{v_j},\,P(\overline{u}_{w_j},\,\overline{t}))\|}-m
\right|\le\ep m\big\}\big|\ge\\
&&\ge1-6\exp(-2\ep^2c(p,\,q)^2m)\to1\ \ (m\to\infty)\,.
\end{eqnarray*}    
\end{cor}
We prove Theorem~\ref{thm_physPrrea} in the next section. 

\section{Proof of Theorem~\ref{thm_physPrrea}}\label{sec_proof}

We deduce Theorem~\ref{thm_physPrrea}
from the next proposition.

\begin{prop}\label{prop_ranSetHalf}
Suppose that $(P,\overline{S})$, where 
$$
\overline{S}=((\overline{s}_1,\,\overline{u}_1),\,(\overline{s}_1,\,\overline{u}_2),\,
\ds,\,\ds,\,(\overline{s}_M,\,\overline{u}_N))\,,
$$ 
are the two $(p,q,n)$-proofreaders and that $a>0$ is a~real number. Then
for every text $\overline{t}\in B_n$
with $\|\overline{t}\|=m\ge1$ errors the following holds.
\begin{enumerate}
\item We have
$$
{\textstyle
\frac{1}{MN}\cdot\big|\big\{j\in[MN]\cc\;
\big|\|P(\overline{s}_{v_j},\,
\overline{t})\|-pm\big|\le a\big\}\big|\ge1-2\exp(-\frac{2a^2}{m})\,.
}
$$
\item We have
$$
{\textstyle
\frac{1}{MN}\cdot\big|\big\{j\in[MN]\cc\;
\big|\|P(\overline{u}_{w_j},\,
\overline{t})\|-qm\big|\le a\big\}\big|\ge1-2\exp(-\frac{2a^2}{m})\,.
}
$$
\item We have
\begin{eqnarray*}
&&{\textstyle
\frac{1}{MN}\cdot\big|\big\{j\in[MN]\cc\;
\big|\|P(\overline{s}_{v_j},\,P(\overline{u}_{w_j},\,\overline{t}))\|-pqm
\big|\le a\big\}\big|}\\
&&{\textstyle\ge1-2\exp(-\frac{2a^2}{m})\,.
}
\end{eqnarray*}
\end{enumerate}
\end{prop}

\noindent
We prove it at the end of this section. Now we 
deduce Theorem~\ref{thm_physPrrea}.

\medskip\noindent
{\bf Deduction of 
Theorem~\ref{thm_physPrrea} from Proposition~\ref{prop_ranSetHalf}. }Suppose that $(P,\overline{S})$, where 
$$
\overline{S}=((\overline{s}_1,\,\overline{u}_1),\,(\overline{s}_1,\,\overline{u}_2),\,
\ds,\,\ds,\,(\overline{s}_M,\,\overline{u}_N))\,,
$$ 
are two $(p,q,n)$-proofreaders, that $\overline{t}\in 
B_n$ is a~text with 
$\|\overline{t}\|=m\ge1$ errors and that $a\in[0,m]$ is a~real number. We set 
$${\textstyle
b=c(p,\,q)\cdot a=\frac{1}{2}\big(\frac{1}{p}+\frac{1}{q}+\frac{2}{pq}\big)^{-1}\cdot a\ \ (\le\frac{pqm}{2})\,.
}
$$
Let $A$ be the set of $j$ ($\in[MN])$ such that 
$|\,\|P(\overline{s}_{v_j},\overline{t})\|-pm|\le b$, $B$ be the set of 
$j$ such that $|\,\|P(\overline{u}_{w_j},\overline{t})\|-qm|\le b$, $C$ be the set of $j$ such that 
$$
|\,\|P(\overline{s}_{v_j},\,P(\overline{u}_{w_j},\,\overline{t}))\|-pqm|\le b
$$ 
and, finally, let $D$  be the set we are interested in, the set of $j$ such that 
$$
|R-m|:=
\big|\,\|P(\overline{s}_{v_j},\,\overline{t})\|\cdot\|P(\overline{u}_{w_j},\,\overline{t})\|\cdot
\|P(\overline{s}_{v_j},\,P(\overline{u}_{w_j},\,\overline{t}))\|^{-1}-m\big|\le a\,,
$$ 
where $\|\ds\|\cdot\|\ds\|\cdot\|\ds\|^{-1}=+\infty$ if $\|P(\overline{s}_{v_j},P(\overline{u}_{w_j},\overline{t}))\|=0$.
Thus if $j\in A$, resp. $j\in B$, then 
$$
\text{$\|P(\overline{s}_{v_j},\,\overline{t})\|=(1+\de_1)pm$, resp.
$\|P(\overline{u}_{w_j},\,\overline{t})\|=(1+\de_2)qm$}\,,
$$
where $|\de_1|\le\frac{b}{pm}$, resp. $|\de_2|\le\frac{b}{qm}$. Similarly,  if $j\in C$ then
$${\textstyle
\|P(\overline{s}_{v_j},\,P(\overline{u}_{w_j},\,\overline{t}))\|=
(1+\de_3)pqm,\ |\de_3|\le\frac{b}{pqm}\le\frac{1}{2}\,.
}
$$
If $j\in A\cap B\cap C$ then
\begin{eqnarray*}
|R-m|&=&{\textstyle
\big|\frac{(1+\de_1)(1+\de_2)}{1+\de_3}m-m \big|}\\
&=&{\textstyle
m\left|
(1+\de_1+\de_2+\de_1\de_2)\sum_{k=0}^{\infty}(-1)^k\de_3^k-1
\right|}\\
&=&{\textstyle m\left|
\sum_{k=1}^{\infty}(-1)^k\de_3^k+
(\de_1+\de_2+\de_1\de_2)\sum_{k=0}^{\infty}(-1)^k\de_3^k
\right|}\\
&\le&m\left(2|\de_3|+
2(|\de_1|+|\de_2|+|\de_1|\cdot|\de_2|)\right)\\
&\le&2b{\textstyle\big(
\frac{1}{pq}+\frac{1}{p}+\frac{1}{q}+
\frac{1}{pq}\big)=a
}
\end{eqnarray*}
and therefore $D\supset A\cap B\cap C$. By items 1--3 of Proposition~\ref{prop_ranSetHalf},
$${\textstyle
\frac{1}{MN}\cdot|D|\ge1-\frac{|A^c|+|B^c|+|C^c|}{MN}\ge1-
6\exp(-\frac{2b^2}{m})=
1-6\exp(-\frac{2c(p,\,q)^2a^2}{m})
}
$$
(for $X\sus[MN]$ we write $X^c$ for the complement $[MN]\setminus X$).
\eproof

We develop tools for the proof of  Proposition~\ref{prop_ranSetHalf}. Let $(P,\overline{S})$, where 
$$
\overline{S}=((\overline{s}_1,\,\overline{u}_1),\,(\overline{s}_1,\,\overline{u}_2),\,
\ds,\,\ds,\,(\overline{s}_M,\,\overline{u}_N))\,,
$$ 
be two $(p,q,n)$-proofreaders. We consider the finite uniform probability space 
$$
\mathrm{PS}(MN)=([MN],\,
\mathcal{P}([MN]),\,\mathrm{Pr})\,,
$$
where $\mathrm{Pr}(\{j\})=
\frac{1}{MN}$, and for $i\in[n]$ the three $n$-tuples of $\mathrm{PS}(MN)$-events 
$$
A_i=\{j\in[MN]\cc\;(\overline{s}_{v_j})_i=1\},\ 
A_i'=\{j\in[MN]\cc\;(\overline{u}_{w_j})_i=1\}
$$
and
$$
C_i=\{j\in[MN]\cc\;(\overline{s}_{v_j})_i=(\overline{u}_{w_j})_i=1\}\,.
$$
By items 1--3 of Proposition~\ref{prop_twoProof}, each of these three $n$-tuples consists of independent events. Also,  $\mathrm{Pr}(A_i)=p$,
$\mathrm{Pr}(A_i')=q$ and
$\mathrm{Pr}(C_i)=pq$
for every $i\in[n]$. Recall that $B^c$ denotes the 
complementary event to the event $B$. For example,  
$$
A_i^c=[MN]\setminus A_i=
\{j\in[MN]\cc\;(\overline{s}_{v_j})_i=0\}\,.
$$
Let $X\sus[n]$ be a~nonempty set.
For any set $Y\sus X$ we consider the three $\mathrm{PS}(MN)$-events
$$
{\textstyle
A_{Y,X}=\bigcap_{i\in Y}A_i\cap\bigcap_{i\in X\setminus Y}A_i^c,\ 
A_{Y,X}'=\bigcap_{i\in Y}A_i'\cap\bigcap_{i\in X\setminus Y}(A_i')^c
}
$$
and
$$
{\textstyle
C_{Y,X}=\bigcap_{i\in Y}C_i\cap\bigcap_{i\in X\setminus Y}C_i^c\,.
}
$$

We employ the following theorem which follows from \cite[Corollary 
A.1.7 on p.~325]{PM} which is proven in \cite[pp.~321--325]{PM}.

\begin{thm}\label{thm_larDev}
Let $n\in\N$, $X\sus[n]$ with 
$|X|=m\ge1$ and let $a>0$ be a~real number. With the above notation we 
have the following bounds.
\begin{enumerate}
\item
$$
{\textstyle
\sum_{\substack{Y,\,Y\sus X\\|\,|Y|-pm\,|\le a}}\mathrm{Pr}(A_{Y,X})\ge1-2\exp(-\frac{2a^2}{m})\,.
}
$$
\item
$$
{\textstyle
\sum_{\substack{Y,\,Y\sus X\\|\,|Y|-qm\,|\le a}}\mathrm{Pr}(A_{Y,X}')\ge1-2\exp(-\frac{2a^2}{m})\,.
}
$$
\item
$$
{\textstyle
\sum_{\substack{Y,\,Y\sus X\\|\,|Y|-pqm\,|\le a}}\mathrm{Pr}(C_{Y,X})\ge1-2\exp(-\frac{2a^2}{m})\,.
}
$$
\end{enumerate} 
\end{thm}
We explain in detail the deduction of this theorem from \cite[Corollary 
A.1.7]{PM}. Let $m\in\N$, $p\in[0,1]$, 
$X_i$ for $i\in[m]$ be $m$ independent two-valued
real random variables in 
a~probability space, with the common distribution $\mathrm{Pr}(X_i=1-p)=p$ 
and $\mathrm{Pr}(X_i=-p)=1-p$, and let
$$
X_0=X_1+X_2+\ds+X_m\,.
$$
The random variable $X_0$ has the centered binomial distribution 
$B(m,p)-pm$. The estimate \cite[Corollary A.1.7]{PM}, more precisely a~particular case of it with all $p_i=p$, says the following.

\begin{thm}[Corollary A.1.7 in \cite{PM}]\label{thm_corInPM}
With the above notation we have for every real number $a>0$ that
$${\textstyle
\mathrm{Pr}(|X_0|>a)<\exp(-\frac{2a^2}{m})\,.
}
$$
\end{thm}
What does the independence of random variables 
$X_i$, $i\in[m]$, which is needed in Theorem~\ref{thm_corInPM}, precisely 
mean? Some events $C_1$, $C_2$, $\ds$, $C_m$ in a~probability space are independent, if for every set $I\sus[m]$ we have
$$
{\textstyle
\mathrm{Pr}\big(\bigcap_{i\in I}C_i\big)=\prod_{i\in I}\mathrm{Pr}(C_i)\,.
}
$$
These events are {\em strongly independent}, if every $m$-tuple of events $D_1$, $D_2$, $\ds$, $D_m$, where each $D_i\in
\{C_i,C_i^c\}$, is independent. A~well known lemma, whose proof we
leave for the interested reader as an exercise, says that any 
independent events are strongly
independent. Discrete (i.e., with at most countable images) real random
variables $Y_1$, $Y_2$, $\ds$, $Y_m$ are independent if for every 
$m$-tuple 
$$
\overline{a}=(\overline{a}_1,\,\overline{a}_2,\,
\ds,\,\overline{a}_m)\in\R^m\,, 
$$
the events
$\mathrm{Pr}(Y_i=\overline{a}_i)$, $i\in[m]$, 
are independent. If the $Y_i$ are two-valued, the image of $Y_i$ being 
$\{a_i,b_i\}$ with $a_i\ne b_i$, then, due to the 
mentioned lemma,  the random variables $Y_i$ are independent if for any single $m$-tuple $\overline{c}\in\R^m$ with $\overline{c}_i\in
\{a_i,b_i\}$ the events
$\mathrm{Pr}(Y_i=c_i)$, $i\in[m]$, 
are independent.

We deduce Theorem~\ref{thm_larDev} from Theorem~\ref{thm_corInPM}. We 
consider only item~1, the arguments
for items~2 and~3 are
similar. So let $(P,\overline{S})$, with 
$$
\overline{S}=((\overline{s}_1,\,\overline{u}_1),\,(\overline{s}_1,\,\overline{u}_2),\,
\ds,\,\ds,\,(\overline{s}_M,\,\overline{u}_N))\,,
$$ 
be as above two $(p,q,n)$-proofreaders and let $X\sus[n]$ with 
$|X|=m\ge1$. We may
assume that $X=[m]$. For $i\in[m]$
we define two-valued 
$\mathrm{PS}(MN)$-random variables $X_i$ by setting ($j\in[MN]$) 
$$
\text{$X_i(j)=1-p$ if $j\in A_i$, 
and $X_i(j)=-p$ if $j\in[MN]\setminus A_i$}\,. 
$$
We set $X_0=\sum_{i=1}^m X_i$. 
Then for any subset $Y\sus X$ we have
$$
\mathrm{Pr}(A_{Y,X})=
\mathrm{Pr}(X_0=|Y|-pm)\,.
$$
By item~1 of Proposition~\ref{prop_twoProof} and the above discussion of independence, 
the random variables $X_i$, $i\in[m]$, are independent. Let $a>0$ 
be a~real number. Applying Theorem~\ref{thm_corInPM} to the 
random variable $X_0$ we get the estimate in item~1 of 
Theorem~\ref{thm_larDev}. 

We proceed to the proof of Proposition~\ref{prop_ranSetHalf}. 

\medskip\noindent
{\bf Proof of Proposition~\ref{prop_ranSetHalf}. }We consider only item~1, the arguments for items~2 and~3 are similar. Suppose that $(P,\overline{S})$, where 
$$
\overline{S}=((\overline{s}_1,\,\overline{u}_1),\,(\overline{s}_1,\,\overline{u}_2),\,
\ds,\,\ds,\,(\overline{s}_M,\,\overline{u}_N))\,,
$$ are two $(p,q,n)$-proofreaders and that $a>0$ is a~real 
number. Let $\overline{t}\in B_n$ with $\|\overline{t}\|=m\ge1$. We may 
assume that the $1$s of $\overline{t}$
are on the first $m$ positions. 
We consider the $\mathrm{PS}(MN)$-events $A_i$, $i\in[m]$. Let $Y\sus[m]=X$. 
Then 
$$
{\textstyle
\frac{1}{MN}|\{j\in[MN]\cc\;\text{for $i\in[m]$ we have $(\overline{s}_{v_j})_i=1$ iff $i\in Y$}\}|=\mathrm{Pr}(A_{Y,X})\,.
}
$$
Thus for any $k\in\N_0$ we have
$${\textstyle
\frac{1}{MN}\cdot\big|\big
\{j\in[MN]\cc\;
\|P(\overline{s}_{v_j},\,\overline{t})\|=k\big\}\big|
=\sum_{\substack{Y,\,Y\sus X\\|Y|=k}}\mathrm{Pr}(A_{Y,X})\,.
}
$$
The stated estimate therefore follows from item~1 of 
Theorem~\ref{thm_larDev}.
\eproof

\section{Concluding remarks}\label{sec_concl}

An interesting direction for further research is the following.
Suppose that 
$(P,\overline{S})$, where 
$$
\overline{S}=((\overline{s}_1,\,\overline{u}_1),\,(\overline{s}_1,\,\overline{u}_2),\,
\ds,\,\ds,\,(\overline{s}_M,\,\overline{u}_N))\,,
$$
are two $(p,q,n)$-proofreaders and that $\overline{t}\in B_n$ is 
a~text. We do not have access to $\overline{t}$ and to 
$m=\|\overline{t}\|$, but the multiset of triples of values
$$
V=\{(P(\overline{s}_{v_j},\,\overline{t}),\,P(\overline{u}_{w_j},\,\overline{t}),\,P(\overline{s}_{v_j},\,P(\overline{u}_{w_j},\,\overline{t})))\in\N_0^3\cc\;j\in[MN]\}
$$
is known. Can we recover the number of errors $m$, or at least some 
bounds on $m$, from $V$? The parameter $n$ is of course known but 
we get two versions of this problem depending on whether $p$ and $q$ are 
known or not. 

Definitions~\ref{def_pnPrRe} and \ref{def_twoPrRe} with accompanying
Propositions~\ref{prop_runs} and \ref{prop_twoProof}
are probably more important results
than Theorem~\ref{thm_physPrrea} whose deduction from the well known 
Theorem~\ref{thm_corInPM} is relatively straightforward. The two 
definitions provide {\em realistic models} (RM) of P\'olya's 
proofreaders. Indeed, if 
$$
\overline{s}=(\overline{s}_1,\,\overline{s}_2,\,\ds,\,\overline{s}_N)
$$
are runs of a~$(p,n)$-proofreader $(P,\overline{s})$ and
$\overline{t}\in B_n$ is a~text, then we can, at least in principle, 
physically go through all $N=l^n$ binary $n$-tuples 
$\overline{s}_j$ and compute all numbers 
$\|P(\overline{s}_j,\overline{t})\|$ of errors the proofreader detects 
in the text during all $N$ runs.

For this to be possible, it is crucial that $p\in\Q$. If 
$p=1/\sqrt{2}$, say, then no RM for a~proofreader with efficiency $p$ 
exists. For irrational $p$ one has to be content with {\em 
probabilistic models} (PM), of random subsets containing (independently) a~given 
element with probability $p$. Their  
construction is easy, but they cannot be physically realized in the 
same sense in which RM can. 

Some paradoxes in probability theory arise when a~``practical'' problem 
is modeled right away with a~PM, 
while one should start carefully with a~RM. The best 
example of this is the well known Bertrand paradox. This paradox 
begins as a~``practical'' problem: what is the probability that 
a~random chord in the unit circle is longer than the length $\sqrt{3}$ of 
the side of the inscribed equilateral triangle. Solutions 
start with various PMs, which automatically assume dimensionless 
points and chords with zero width. This is a~mistake, no such objects 
exist in reality. Instead, one should start with a~RM where points 
have size $s>0$ and chords have width $w>0$. Once such RM is 
established, one can send $s,w\to0$ and see what happens in the limit. This is the 
approach in \cite{klaz}, which we hope to rework and improve in the 
near future.

\bigskip\noindent
{\em Department of Applied Mathematics\\ 
Faculty of Mathematics and Physics\\
Charles University\\ 
Malostransk\'e n\'am\v est\'\i\ 25\\
118 00 Praha\\
Czechia}\\

\end{document}